\documentclass[12pt,reqno]{amsart}

\usepackage{amsfonts,color,amsthm,amsmath,amssymb}

\usepackage{color}
\def\Box{\vcenter{\vbox{\hrule\hbox{\vrule
     \vbox to 8.8pt{\hbox to 10pt{}\vfill}\vrule}\hrule}}}

\def\qed{{\hfill$\square$}}
\def\proof{{\vspace{-0.0cm}\bf Proof: \,}}
\def\N{{\mathbb N}}

\def\F{{\mathbb F}}

\def\mod{{\mathrm{mod\,\,}}}

\def\max{{\mathrm{max}}}
\def\min{{\mathrm{min}}}

\newtheorem{theorem}{Theorem}[section]

\newtheorem{lemma}[theorem]{Lemma}
\newtheorem{remark}[theorem]{Remark}

\newtheorem{corollary}[theorem]{Corollary}

\newtheorem{proposition}[theorem]{Proposition}

\numberwithin{equation}{section}

\begin{document}
\title[Upper bounds on the size of transitive subtournaments in digraphs]
{Upper bounds on the size of transitive subtournaments in digraphs}

\author[Koji Momihara and Sho Suda]{Koji Momihara$^*$ and  Sho Suda$^{\dagger}$}

\thanks{$^\ast$ Koji Momihara is supported by JSPS KAKENHI Grant Number (C)24540013.}
\thanks{$^{\dagger}$
Sho Suda is supported by  JSPS KAKENHI Grant Number 15K21075.}

\address{$^\ast$ Faculty of Education, Kumamoto University, 2-40-1 Kurokami, Kumamoto 860-8555, Japan} \email{momihara@educ.kumamoto-u.ac.jp}

\address{$^{\dagger}$ Department of Mathematics Education, Aichi University of 
Education, 1 Hirosawa, Igaya-cho, Kariya, Aichi 448-8542, Japan} \email{suda@auecc.aichi-edu.ac.jp}


\begin{abstract}
In this paper, we consider upper bounds on the size of transitive subtournaments 
in a digraph. 
In particular, we give an analogy of Hoffman's bound for the size of cocliques in a regular graph. 
Furthermore, we partially improve the Hoffman type bound for doubly regular 
tournaments by using the technique of Greaves and Soicher for strongly regular graphs~\cite{GS}, which gives a new application of  block intersection polynomials.  
\end{abstract}

\maketitle

\section{Introduction}
The problem to find a sharp bound on the size of cliques or cocliques in a graph has been well-studied, and some nontrivial bounds  have been known based on linear algebraic techniques, cf.~\cite{Ha1,Ha2}. In particular, it is well-known that  if a $k$-regular graph with $v$ vertices has a coclique 
of size $s$, then  
\begin{equation}\label{eq:hoff:ori}
s\le -v\lambda_{\min}/(k-\lambda_{\min}),
\end{equation}
where $\lambda_{\min}$ is the 
minimum eigenvalue of the adjacency matrix of $G$. 
This bound is an unpublished result of A. J. Hoffman, and so this is known as {\it Hoffman's bound}. 
Furthermore, the case where 
the equality holds in this bound was studied in relation to some combinatorial structures, such as  strongly regular graphs and association schemes~\cite{Ha2}.  
The Hoffman bound for strongly regular graphs is sometimes referred to as 
the Delsarte bound. 
The Delsarte bound is difficult to improve in general. In fact, until recently, there was only a few classes of strongly regular graphs, for which the bound is improvable.  
For example, Bachoc, Matolcsi and Ruzsa~\cite{BMR} improved the Delsarte bound for Paley graphs of nonsquare 
order using the properties of quadratic residues of finite fields. 
Very recently, it was announced that Greaves and Soicher~\cite{GS} improved the Delsarte bound for a large class of strongly regular graphs using block-intersection polynomials. 

In this paper, we are interested in the size of transitive subtournaments in a 
digraph. The problem to find a sharp  bound on the size of transitive subtournaments in a tournament was initially considered by Erd\H{o}s and Moser~\cite{EM} in 1964.  It is clear that any tournament with $v$ vertices contains a transitive 
subtournament with $1+\lfloor \log_2(v)\rfloor$ vertices. 
This bound is tight for the case where $v=7$. In fact, the Paley 
tournament on seven vertices contains no transitive tournament with four vertices. They also proved that there are tournaments with $v$ vertices without a transitive tournament of size 
$2(1+\lfloor \log_2(v)\rfloor)$.  On the other hand, Reid and Parker~\cite{RP} showed that the lower bound above is not tight for large $v$.  There have been some studies of upper and lower bounds on the size of transitive subtournaments in a digraph to improve Erd\H{o}s-Moser's bound,  
cf. \cite{Moo,Sa1,Sa2,T}. 
In particular, in \cite{T}, it was shown by a simple counting argument that for the maximum size $s$ of transitive subtournaments in 
a doubly regular tournament with $v$ vertices 
\begin{equation}\label{eq:bound}
s\le  \frac{-3+\sqrt{13+12v}}{2}. 
\end{equation}

In this paper, we consider analogies of  
the Hoffman bound and  the Greaves-Soicher bound to digraphs. 
In particular, as an analogy of the Hoffman bound for regular graphs, we show that if a regular digraph $G$ with $v$ vertices has a transitive subtournament of size $s$, then 
\[
s\le 
 \frac{-3\theta_{\max}^2+\sqrt{9\theta_{\max}^4+4v^2+12\theta_{\max}^2v^2}}{2v}, 
\] 
where $\theta_{\max}$ is the 
maximum eigenvalue of the Seidel matrix of $G$. As an immediate corollary, we have for a 
doubly regular tournament $
s \le (-3+\sqrt{13+12v})/2$, which coincides with the bound~\eqref{eq:bound}. Furthermore,  this bound can be partially improved by using 
Greaves-Soicher's technique for strongly regular graphs~\cite{GS}, which gives  a new application of  block-intersection polynomials.   

This paper is organized as follows. 
In Section~\ref{sec:2}, we introduce basic facts on spectra of digraphs. 
In Section~\ref{sec:3}, we obtain an upper bound on the size of transitive subtournaments in a  
digraph using the ``interlacing property'' of eigenavalues without any restriction. 
Section~\ref{sec:hof} is  the main part of this paper, where we consider an analogy of the Hoffman bound to 
digraphs. 
In Section~\ref{sec:5}, we introduce the adjacency polynomials for doubly regular tournaments, and partially improve
the Hoffman type bound for doubly regular tournaments, that is, the bound~\eqref{eq:bound}. Finally, we give some open problems related to this 
study in Section~\ref{sec:6}. 
\section{Preliminaries}\label{sec:2}
Let $G=(V,E)$ be a digraph with $v$ vertices.  An {\it adjacency matrix} $A$ of $G$ is a $v\times v$ matrix whose columns and rows are labeled by the vertices of $G$, and its entries are defined by 
\begin{align*}
A_{x,y}=\begin{cases}
1,\quad &\textup{if $(x,y)\in E$},\\
0,\quad & \textup{otherwise}.
\end{cases}
\end{align*}
It is clear that $A+A^\top $ is symmetric and $\sqrt{-1}(A-A^\top )$ is Hermitian, where $A^\top$ denotes the transpose of $A$. 
The matrix $S_{G}:=\sqrt{-1}(A-A^\top )$ is called the {\it Seidel matrix} of $G$. 
Let ${\bf 1}$ be the all-one vector of length $v$. 
If $|\{x\in V:(x,y)\in E\}|=|\{x\in V:(y,x)\in E\}|=k$ for any $y\in V$, then  the digraph is called {\it regular}. This condition is equivalent to the fact that $A\cdot {\bf 1}=A^\top \cdot {\bf 1}$ and all entries of the vector $A\cdot {\bf 1}$ are equal to $k$. 


A digraph is called a {\it tournament} if $A+A^\top =J-I$, where $J$ is the all-one square matrix of order $v$ and $I$ is the identity matrix of order $v$. A tournament $G$ is called  {\it doubly regular} if it is $(v-1)/2$-regular  and the number of vertices dominated by a pair of two distinct vertices simultaneously is constant, say $\lambda$, not depending on the choice of the pair. It is clear that $\lambda=(v-3)/4$, which implies that $v\equiv 3\,(\mod{4})$. The adjacency matrix $A$ satisfies that 
\[
AA^\top =\frac{v+1}{4}I+\frac{v-3}{4}J. 
\]

Noting that $A+A^\top =J-I$, $A$ has   the eigenvalues $\frac{v-1}{2},\frac{-1+\sqrt{-v}}{2},\frac{-1-\sqrt{-v}}{2}$. In other words, $S_G$ has 
$0,\sqrt{v},-\sqrt{v}$ as its eigenvalues \cite[Theorem~2.5]{NS}.

Since $S_G$ is Hermitian, $S_G$ has only real eigenvalues. Let $E_\theta$ be 
the orthogonal projection matrix of an eigenvalue $\theta$. Then, 
$S_G$ has  the spectral decomposition 
\begin{equation}\label{eq:spe}
S_G=\sum_{\theta\in \mbox{ev}(G)}\theta E_\theta, 
\end{equation}
where $\mbox{ev}(G)$ is the set of distinct eigenvalues of $S_G$. 
We denote by $m_\theta$ the multiplicity corresponding to $\theta\in\mbox{ev}(G)$. 
The {\it main angle} $\beta_\theta$ is defined to be $\beta_\theta=(1/\sqrt{v})||E_\theta\mathbf{1}||^2$. 
A {\it main eigenvalue}  of $S_G$ is an eigenvalue $\theta$ with $\beta_\theta\neq0$. 

Since $S_G$ is Hermitian and $S_G^\top =-S_G$, $E_\theta$'s satisfy the following basic properties.  
\begin{enumerate}
\item $E_\theta^2=E_{\theta}$ and $E_\theta E_{\tau}=O$ for $\theta \not=\tau$,
\item $S_G E_\theta=\theta E_\theta$,
\item $\sum_{\theta \in\mbox{ev}(G)}E_\theta=I$, 
\item if $\theta\in\mbox{ev}(G)$, then $-\theta\in\mbox{ev}(G)$, 
\end{enumerate}
where $O$ denotes the zero matrix of order $v$. 
 
Define the set  $M(G)$ of main eigenvalues to be $M(G)=\{\theta\in \mathrm{ev}(G): \beta_\theta\neq 0\}$, and the matrix $F_G$ to be $F_G=\sum_{\theta\in M(G)}E_\theta$. 
Note that $M(G)\neq\emptyset$ since $\sum_{\theta\in\mathrm{ev}(G)}\beta_{\theta}^2=1$.

In the rest of this paper, we are interested in the size of transitive subtournaments in a digraph. A transitive tournament is a tournament satisfying the following: if $(x,y)\in E$ and $(y,z)\in E$, then $(x,z)\in E$. 
After reordering the vertices appropriately, we may assume that a 
transitive tournament has the adjacency matrix
\begin{align}\label{eq:SeidelTra}
A=\begin{pmatrix}
0 & 0&0&\cdots&0 \\
1 & 0&0&\cdots&0\\
1&1&0&\cdots&0\\
\vdots & \vdots&\vdots &\ddots&\vdots\\
1&1&1&\cdots&0
\end{pmatrix}.
\end{align} 

\section{Interlacing}\label{sec:3}
In this section, we will use the following well-known fact on interlacing of eigenvalues. 
\begin{proposition}{\em (\cite[Theorem 4.3.17]{HJ})}
Let  $A$ be a Hermitian matrix of order $v$ with eigenvalues 
\[
\lambda_1\ge \lambda_2\ge \cdots\ge \lambda_v. 
\]
Let $B$ be a principal submatrix of $A$ of  order $m$ 
with eigenvalues 
\[
\mu_1\ge \mu_2\ge \cdots\ge \mu_m. 
\]
Then, the eigenvalues of $B$ interlace those of $A$, i.e., 
\[
\lambda_i\ge \mu_i\ge \lambda_{v-m+i}, \, \, \, \, i=1,2,\ldots,m.
\]
\end{proposition}

Using interlacing of eigenvalues, we obtain the following theorem. 
\begin{theorem}\label{thm;inte}
Let $G$ be a digraph with $v$ vertices and  $\Gamma$ be a 
transitive subtournament of size $s$ in $G$. Let $\theta_i$, 
$i=1,2,\ldots,v$, be the eigenvalues of $S_G$ with ordering $\theta_1\geq\cdots\geq\theta_v$. Then, 
\begin{equation}\label{eq:inter}
s\leq \frac{(2i-1)\pi}{2\mathrm{arccot}(\theta_{i})}
\end{equation}
for $i=1,\ldots,\lfloor s/2 \rfloor$. 
\end{theorem}
\proof
It is easily shown that the eigenvalues of $S_\Gamma$ are $\cot(\frac{(2i-1)\pi }{2s})$, $i=1,\ldots,s$.
By interlacing for $S_G$ and $S_\Gamma$, we have $\cot(\frac{(2i-1)\pi}{2s})\leq \theta_{i}$, that is, the inequality \eqref{eq:inter}. 
\qed

\vspace{0.3cm}
By applying Theorem~\ref{thm;inte} to doubly regular tournaments, we have the following corollary. 
\begin{corollary}\label{cor;inte}
Let $G$ be a doubly regular tournament with $v$ vertices and  $\Gamma$ be  a transitive subtournament of size $s$ in $G$. Then, 
\[
s\leq \frac{\pi}{2\mathrm{arccot}(\sqrt{v})}. 
\]
\end{corollary}
\proof This corollary follows by the fact that $\theta_i=\sqrt{v}$ for $i=1,\ldots,(v-1)/2$. 
\qed

\vspace{0.3cm}
We list the values of $\frac{\pi}{2\mathrm{arccot}(\sqrt{v})}$ for small $v\equiv 3\,(\mod{4})$ in Table~\ref{tableEx3}. As far as we checked for small cases, this bound is not better than the bound obtained in Section~\ref{sec:hof}.  On the other hand, 
the advantage of Theorem~\ref{thm;inte} is that the result is applicable 
to general digraphs not computing $\beta_\theta$'s.  
\begin{table}[h]
\begin{center}
\begin{tabular}{|c||c|c|c|c|c|c|c|c|c|c|}
\hline 
$v$ &   7&11&15&19&23&27&31&35\\
\hline 
$\frac{\pi}{2\mathrm{arccot}(\sqrt{v})}$&   4.346&5.363&6.216&6.965&7.641&8.261&8.839&9.380 \\
\hline 
\end{tabular}
\end{center}
\caption{The values of $\frac{\pi}{2\mathrm{arccot}(\sqrt{v})}$ for small $v$.}\label{tableEx3}
\end{table}
\section{Analogy of Hoffman's bound}\label{sec:hof}
In this section we consider an analogy of Hoffman's bound to digraphs. 
Recall that $M(G):=\{\theta\in \mathrm{ev}(G): \beta_\theta\neq 0\}$ and $F_G:=\sum_{\theta\in M(G)}E_\theta$. 
\begin{lemma}\label{lem:chi1}
Let $G$ be a digraph with $v$ vertices 
and $\chi$ be a $(0,1)$-vector of norm $s$. 
Then, $\chi^\top F_G\chi\geq\frac{1}{v}s^2$. 
\end{lemma}
\proof
By the definition of main angles, the space $F_G\mathbb{R}^v$ contains $\text{span}\{{\bf 1}\}$. 
Letting $E$ be the orthogonal projection onto $F_G\mathbb{R}^v\cap \text{span}\{{\bf 1}\}^\perp$,  
we have
\begin{align*}
F_G=\frac{1}{v}J+E,
\end{align*}
from which the result follows. 
\qed
\begin{lemma}\label{lem:chi2}
Let $G$ be a digraph with $v$ vertices and  $\Gamma$ be a transitive subtournament of size $s$ in $G$. 
Then it holds that
\[
s(s^2-1)/3\leq \chi^\top S_GS_G^* \chi, 
\]
with equality if and only if for any vertex $x\not\in V(\Gamma)$, the number of vertices in $V(\Gamma)$ dominating $x$ equals to the number of vertices in $V(\Gamma)$ dominated by $x$. 
\end{lemma}
\proof
After reordering the vertices appropriately, we may set 
\begin{align*}
S_G=\begin{pmatrix}
S_\Gamma & S_{12} \\
S_{21} & S_{22}
\end{pmatrix} .
\end{align*} 
Then it is easy to see that 
\begin{align*}
0\le {\bf 1}^\top S_{12} S_{12}^\ast{\bf 1}= \chi^\top S_G S_G^\ast \chi-{\bf 1}^\top S_ \Gamma S_\Gamma^\ast {\bf 1}=\chi^\top S_G S_G^\ast \chi-s(s^2-1)/3.  
\end{align*}
The equality holds if and only if $S_{12}^*{\bf 1}={\bf 0}$, which is equivalent to the desired condition.   
\qed

\vspace{0.3cm}
The following is our main theorem in this section. 
\begin{theorem}\label{thm:hoff1}
Let $G$ be a digraph with $v$ vertices and   $\Gamma$ be a  transitive subtournament of size $s$ in $G$.  
Let $\alpha=\max\{\theta:\theta\in M(G)\}$ and  $\gamma=\max \{\theta:\theta\in\mathrm{ev}(G)\setminus M(G)\}$.
If $\alpha\leq \gamma$, then 
\begin{align}\label{ineq:1}
 s\le \frac{3\alpha^2-3 \gamma^2 + \sqrt{4v^2(1+3\gamma^2)+9(\gamma^2-\alpha^2)^2}}{2 v}.
\end{align}
\end{theorem}
\proof
Let $\chi$ be the characteristic vector of $\Gamma$. 
Then 
\begin{align}
\chi^\top S_G S_G^\ast \chi=&\, \chi^\top \big( \sum_{\theta\in\mbox{ev}(G)}\theta E_\theta\big)
\big( \sum_{\theta\in\mbox{ev}(G)}\theta E_\theta^\ast\big)\chi\nonumber\displaybreak[0]\\
=&\, \chi^\top \big( \sum_{\theta\in\mbox{ev}(G)}\theta^2 E_\theta\big)\chi\nonumber\displaybreak[0]\\
\leq&\, \alpha^2\chi^\top \big( \sum_{\theta\in M(G)}E_\theta\big)\chi+\gamma^2\chi^\top \big( \sum_{\theta\in\mbox{ev}(G)\setminus M(G)} E_\theta\big)\chi\nonumber\displaybreak[0]\\
= &\, \alpha^2 \chi^\top F_G\chi+\gamma^2\chi^\top \big( I-F_G )\chi \nonumber\displaybreak[0]\\
=&\, (\alpha^2-\gamma^2)\chi^\top F_G\chi+\gamma^2\chi^\top \chi \nonumber\displaybreak[0] \\
=&\, (\alpha^2-\gamma^2)\chi^\top F_G\chi+\gamma^2 s. \label{eq:hoff2i}
\end{align}

By the inequality~\eqref{eq:hoff2i} and Lemma~\ref{lem:chi2}, we have 
\begin{align}\label{eq:hoff11}
0\leq (\alpha^2-\gamma^2)\chi^\top F_G\chi+\gamma^2 s-s(s^2-1)/3.
\end{align}

Since $\alpha\leq \gamma$, the desired inequality follows from Lemma~\ref{lem:chi1} and the inequality~\eqref{eq:hoff11}.  
\qed

\vspace{0.3cm}
We now consider the case where $G$ is a regular digraph.   
Then $M(G)=\{0\}$ and $\gamma=\max\{\theta:\theta\in \mathrm{ev}(G)\setminus M(G)\}$ is equal to the maximum eigenvalue of $S_G$. 
Thus we have the following corollary.  
\begin{corollary}\label{cor:regu}
Let $G$ be a regular digraph with $v$ vertices and $\Gamma$ be a transitive subtournament of size $s$ in $G$. 
Then, 
\begin{equation}\label{eq:regu}
s\le  \frac{-3\theta_{\max}^2+\sqrt{9\theta_{\max}^4+4v^2+12\theta_{\max}^2v^2}}{2v}, 
\end{equation} 
where $\theta_{\max}$ is the 
maximum eigenvalue of the Seidel matrix of $G$. 
\end{corollary}
  
Next, we move on to the case where $G$ is a regular tournament.  
Since $A=\frac{-\sqrt{-1}S_G+J-I}{2}$
and  $S_G+S_G^\top =O$, we have 
\[
(A+{A}^\top )\chi=(J-I)\chi=s {\bf 1}-\chi. 
\]
Hence, for any $x\not\in V(\Gamma)$ \[
|\{y\in V(\Gamma): (x,y)\in E\}|+|\{y\in V(\Gamma): (y,x)\in E\}|=s. 
\]
On the other hand, if the equality holds in the inequality~\eqref{eq:regu}, we have 
\[
|\{y\in V(\Gamma): (x,y)\in E\}|=|\{y\in V(\Gamma): (y,x)\in E\}|. 
\]
Therefore, $s$ must be even. 
Summing up, we have the following corollary. 
\begin{corollary}\label{cor:par}
Let $G$ be a regular  
tournament with $v$ vertices. If the equality holds in the inequality~\eqref{eq:regu}, then 
$s$ is even. 
\end{corollary}
Also, by applying Corollary~\ref{cor:regu} to doubly regular tournaments, 
we have the following.  
\begin{corollary}\label{cor:doub}
Let $G$ be a doubly regular  
tournament with $v$ vertices  and  $\Gamma$ be a  transitive subtournament of size $s$ in $G$.  
Then, it holds that 
\begin{align}\label{ineq:drt}
s\le \frac{-3 + \sqrt{13 + 12 v}}{2}. 
\end{align}
\end{corollary}
\proof
The corollary follows from the fact that $\mbox{ev}(S_G)=\{0,\pm \sqrt{v}\}$ \cite[Theorem~2.5]{NS}. 
\qed 

\vspace{0.3cm}
In Table~\ref{tableEx2}, we computed the maximum size of 
transitive subtournaments in a known doubly regular tournament with $v$ vertices for small $v$. In the row of ``maximum'', we found a 
doubly regular tournament with a transitive subtournament having the indicated number of vertices in known doubly regular tournaments given in \cite{Ma}. Here, the upper bound $s\le 6$ for $v=23$ is obtained by Corollary~\ref{cor:par}. 
Also, we list the maximum size of 
transitive subtournaments in the Paley tournament. Here, the Paley tournament  is defined as follows: Let $q$ be a prime power congruent to $3$ modulo $4$, and let $C_0$ be the set of nonzero squares of the finite field $\F_q$. The {\it Paley tournament of order $q$} is the tournament  with the elements of $\F_q$ as vertices; $(x,y)\in E$ 
if and only if $x-y\in C_0$. See \cite{Sa1} for further computational results 
on the maximum size of 
transitive subtournaments in the Paley tournament.  
\begin{table}[h]
\begin{center}
\begin{tabular}{|c||c|c|c|c|c|c|c|c|c|c|}
\hline 
$v$&   7&11&15&19&23&27&31&35 \\
\hline 
$\#$&   1&1&2&2&37&722&$\ge$ 5&$\ge$ 486 \\
\hline 
\hline 
upper bound &  3&4&5&6&6&7&8&8\\
\hline 
maximum &  3&4&5&5&6&7&$\ge$ 7&$\ge$ 7\\
\hline
Paley &  3&4&-&5&5&5&7&-\\
\hline
\end{tabular}
\end{center}
\caption{The size of transitive subtournaments in a doubly regular tournament}\label{tableEx2}
\end{table}

Let $G$ be a doubly regular tournament with $v$ vertex and $x$ be a vertex of $G$. 
Let $G'$ be the induced subgraph with vertex set $V(G)\setminus \{x\}$ of $G$. 
Then the tournament $G'$ has $v-1$ vertices and is biregular.  
Furthermore, by \cite[Theorem~1.1]{NS}, the spectrum of the tournament $G'$ is 
\begin{align*}
\mbox{ev}(G')=\{\pm\sqrt{v},\pm1\},\quad \beta_{\pm\sqrt{v}}=0,\quad \beta_{\pm1}=1/\sqrt{2}. 
\end{align*}
By applying Theorem~\ref{thm:hoff1} as $\alpha=1$ and $\gamma=\sqrt{v}$, we have 
\[
s\le \frac{-3+\sqrt{25+12v}}{2}=\frac{-3+\sqrt{13+12(v+1)}}{2}, 
\]
which coincides with the bound~\eqref{ineq:drt}. 

\section{Bounds from block intersection polynomials}\label{sec:5}
The concept of {\it block intersection polynomials} was first introduced by Cameron and Soicher in \cite{CS}. 

For a non-negative integer $k$, define the polynomial 
\[
P(x,k):=x(x-1)\cdots (x-k+1). 
\]
Thus, for $n$ a non-negative integer, ${n\choose k}=P(n,k)k!$.  
For real number sequences $M=[m_0,\ldots,m_s]$, $\Lambda=[\lambda_0,\ldots,\lambda_t]$ with $t\le s$, define the {\it block intersection polynomial} 
\begin{equation}\label{def:bip}
B(x,M,\Lambda)=\sum_{j=0}^t{t \choose j}P(-x,t-j)\Big(
P(s,j)\lambda_j-\sum_{i=j}^sP(i,j)m_i
\Big). 
\end{equation}
Cameron-Soicher~\cite{CS} and Soicher~\cite{S1} proved the following result. 
\begin{theorem}\label{thm:bip}
Let $s$ and $t$ be non-negative integers with $s\ge t$, let $n_0,\ldots,n_s$, 
$m_0,\ldots,m_s$, and $\lambda_0,\ldots,\lambda_t$ be real numbers, such that 
\[
\sum_{i=0}^{s}{i\choose j}n_i={s\choose j}\lambda_j, \, \, \, j=0,1,\ldots,t, 
\] 
and let $B(x)$ be the block intersection polynomial defined in \eqref{def:bip}. If $m_i\le n_i$ for 
all $i$, then $B(b)\ge 0$ for every integer $b$. 
\end{theorem}
See \cite{S1} for further properties of block intersection polynomials. 

Cameron and Soicher~\cite{CS} discussed the multiplicity of a block in  a 
$t$-design using block intersection polynomials.   Soicher~\cite{S1,S2} defined 
adjacency polynomials for edge-regular graphs as a special form of  block intersection polynomials and discussed the existence of cliques in edge-regular graphs. Very recently, Greaves and Soicher~\cite{GS} announced that 
they improved the Hoffman bound for strongly regular graphs using 
adjacency polynomials. In this section, we define adjacency polynomials for digraphs and give an improved bound of the Hoffman type bound given in Corollary~\ref{cor:doub} for doubly regular tournaments.  
\subsection{Adjacency polynomials for doubly regular tournaments}
Let $G$ be a doubly regular tournament with $v$ vertices. 
Note that $v=4m-1$ for some $m\in \N$.  Then the number of vertices dominated by each vertex is 
$k=2m-1$, and the number of vertices dominated by two distinct  vertices simultaneously is 
$\lambda=m-1$. 

Let $\Gamma$ be a transitive subtournament of $G$ with 
$s$ vertices. For $D\subseteq V(G)$, 
let 
\[
\lambda_D:=|\{q\in V(G)\setminus V(\Gamma):D\subseteq N(q)\}|, 
\]
where $N(q)$ is the set of vertices dominating $q$, 
and for $0\le j\le s$, set 
\[
\lambda_j:=\frac{\sum_{D\in {V(\Gamma)\choose j}}\lambda_D}{{s \choose j}}. 
\]
Then, it is clear that $\lambda_0=|V(G)\setminus V(\Gamma)|=v-s$ and 
\[
\lambda_1=\frac{sk}{s}-\frac{1}{s}\sum_{a\in V(\Gamma)}\mbox{outdeg}_{ \Gamma}(a)=k-\frac{s-1}{2}, 
\]
where $\mbox{outdeg}_{\Gamma}(a)$ is the outdegree of $a\in V(\Gamma)$ in $\Gamma$. 
Furthermore, we have 
\[
\lambda_2=\lambda-\frac{2}{s(s-1)}\sum_{a,b\in V(\Gamma);a\not=b}\mbox{dom}_{\Gamma}(a,b), 
\]
where $\mbox{dom}_{\Gamma}(a,b)$ is the number of vertices of $\Gamma$ dominated by $a$ and $b$. 
Since $\Gamma$ is transitive, we have $\sum_{a,b\in V(\Gamma);a\not=b}\mbox{dom}_{\Gamma}(a,b)={s\choose 3}$. 
Hence, $\lambda_2=\lambda-\frac{s-2}{3}$. 

\begin{lemma}
Let $n_i=|\{q\in V(G)\setminus V(\Gamma):|N(q)\cap V(\Gamma)|=i\}|$, $0\le i\le s$. Then,  for $j=0,1,2$ 
\[
\sum_{i=0}^s{i\choose j}n_j={s\choose j}\lambda_j.
\]
\end{lemma}
\proof
For $j=0,1,\ldots,s$, we count in two ways the number $N_j$ of ordered pairs 
$(q,D)$, where $q\in V(G)\setminus V(\Gamma)$ and $D$ is a $j$-subset of $N(q)\cap V(\Gamma)$. 

Each $j$-subset $D$ of $V(\Gamma)$ contributes exactly $\lambda_D$ pairs of the form $(-,D)$ to $N_j$. Hence, by the definition of $\lambda_j$, 
\[
N_j=\sum_{D\in {V(\Gamma) \choose j}}\lambda_D={s \choose j}\lambda_j. 
\]
On the other hand, each $q\in V(G)\setminus V(\Gamma)$ contributes exactly ${|N(q)\cap V(\Gamma)|\choose  j}$ pairs of the form $(q,-)$ to $N_j$. Hence, by the  definition of $n_i$, 
\[
N_j=\sum_{q\in V(G)\setminus V(\Gamma)}{|N(q)\cap V(\Gamma)|\choose  j}=\sum_{i=0}^s{i \choose j}n_i. 
\]
This completes the proof. 
\qed

\vspace{0.3cm}
By the lemma above, the integers $n_0,\ldots,n_s$, $\lambda_0,\lambda_1,\lambda_2$ satisfy  the assumption of Theorem~\ref{thm:bip}. Define 
\begin{align*}
B(x)=&\,B(x,[0^{s+1}],[v-s,k-\frac{s-1}{2},\lambda-\frac{s-2}{3}])\\
=&\,\sum_{j=0}^2{2\choose j}P(-x,2-j)P(s,j)\lambda_j\\
=&\,x(x+1)(v-s)-2xs(k-(s-1)/2)+s(s-1)(\lambda-(s-2)/3). 
\end{align*}
We now define the {\it adjacency polynomial} for doubly regular tournaments 
with $v=4m-1$ vertices:  
\[
C(x,y)=x(x+1)(4m-1-y)-2xy(2m-(y+1)/2)+y(y-1)(m-(y+1)/3).  
\]
Then, we have the following result by Theorem~\ref{thm:bip}: 
\begin{proposition}\label{thm:bip2}
Let $G$ be a doubly regular tournament with $v=4m-1$ vertices. 
If $G$ contains a transitive  subtournament 
of size $s$, then $C(b,s)\ge 0$ 
for all integers $b$. 
\end{proposition}
By this proposition, if $C(b,s)<0$ for some integer $b$, then $G$ 
can not contain a transitive subtournament 
of size $s$. For example, 
put $y=3$ and $m=2$, and then $C(x,3)=4(x-1)^2\ge 0$. 
On the other hand, we have $C(x,4)=3x^2-9x+4$ and 
$C(1,4)<0$. Hence, a doubly regular tournament with seven vertices can not contain a transitive subtournament 
of size $4$. 
\subsection{Improved bound for doubly regular tournaments}
By solving the equation $C(x,y)=0$ for $x$, we have 
\[
x=\frac{3 (-1 + 4 m - y) (-1 + y) \pm \sqrt{-3 (-1 + 4 m - y) (-1 + y) (-3 + 12 m - 4 y - y^2)}}{6 (-1 + 4 m - y)}. 
\]
Since we can assume that $1<y<4m-1$, $C(x,y)$ is nonnegative  if
$3-12m+4y+y^2\le 0$, i.e., $-2-\sqrt{1+12m}\le y\le -2+\sqrt{1+12m}$. 

Put 
\[
z_{m,y}:=\frac{\sqrt{-3 (-1 + 4 m - y) (-1 + y) (-3 + 12 m - 4 y - y^2)}}{6 (-1 + 4 m - y)}.
\] 
To find a bound on the size of transitive subtournaments, 
we need to find an integral $x$ in the interval $((-1 + y)/2 - z_{m,y},(-1+y)/2 + z_{m,y})$
for some integer $y> -2+\sqrt{1+12m}$. 

\begin{theorem} \label{thm:bip3}
Let $G$ be a doubly regular tournament with $v=4m-1$ vertices and let 
$\Gamma$ be a transitive subtournament of size $s$ in  $G$. 
Let $\epsilon:=\sqrt{1+12m}-\lfloor \sqrt{1+12m}\rfloor$. The, the following hold. 
\begin{itemize}
\item[(1)] $s\le -1+\lfloor \sqrt{1+12m}\rfloor$.  
\item[(2)] 
$s\le -2+\sqrt{1+12m}$ if $\epsilon=0$. 
\item[(3)] $s\le -2+\lfloor\sqrt{1+12m}\rfloor$ if  
$\epsilon\not=0$ and $-1+\lfloor \sqrt{1+12m}\rfloor$ is odd. 
\item[(4)] $s\le -2+\lfloor\sqrt{1+12m}\rfloor$ if 
$\epsilon\not=0$ and  $\lfloor \sqrt{1+12m}\rfloor >    (-1+\sqrt{1+48m})/2$. 
\end{itemize}
\end{theorem}
\proof 
(1) If $\lfloor\sqrt{1+12m}\rfloor$ is even, put $(x,y)=((\lfloor\sqrt{1+12m}\rfloor)/2,\lfloor\sqrt{1+12m}\rfloor)$. Then, 
\[
C(x,y)=-\frac{1}{12} (-\epsilon + \sqrt{1 + 12 m})
(3 - 3 \epsilon + \epsilon^2 + 3 \sqrt{1 + 12 m} - 
   2 \epsilon \sqrt{1 + 12 m})<0.
\]   
If $\lfloor\sqrt{1+12m}\rfloor$ is odd, 
put $(x,y)=((\lfloor\sqrt{1+12m}\rfloor-1)/2,\lfloor\sqrt{1+12m}\rfloor)$. Then, 
\[
C(x,y)=\frac{1}{12} 
(-2 + \epsilon) (-1 - \epsilon + \sqrt{1 + 12 m}) (2 - \epsilon + 
2 \sqrt{1 + 12 m})<0.
\]  
Then, by Proposition~\ref{thm:bip2}, the assertion follows. 

(2) Since $C((-1+\sqrt{1+12m})/2,-1+\sqrt{1+12m})=-m<0$, 
 the assertion follows.

(3) Put $y=-1+\lfloor\sqrt{1+12m}\rfloor$. Since $-1+\lfloor \sqrt{1+12m}\rfloor>-2+\sqrt{1+12m}$ and $(-1+y)/2$ is integer, we have  $C((-1+y)/2,-1+\lfloor \sqrt{1+12m}\rfloor)<0$.

(4) If $2z_{m,y}>  1$, then we can find an integral $x$  in the interval $((-1 + y)/2 - z_{m,y},(-1+y)/2 + z_{m,y})$. 
The bound $2z_{m,y}> 1$ is equivalent to $2-12m+3y+y^2> 0$, 
i.e., 
$y> (-3+\sqrt{1+48m})/2$. 
By substituting $y=-1+\lfloor\sqrt{1+12m}\rfloor$ into this inequality, 
we have $\lfloor \sqrt{1+12m}\rfloor >  (-1+\sqrt{1+48m})/2$. 
The proof is now complete.  
\qed

\vspace{0.3cm}
For example, for $v=27$, we have $\sqrt{1+12m}=9.21954...$. 
In this case, the condition of Theorem~\ref{thm:bip3} (4) is satisfied, and we have $s\le 7$. 
\begin{remark}
Theorem~\ref{thm:bip3} (3) is particularly important because 
the result can improve the bound~\eqref{ineq:drt}. For example, 
in the case where $v=71$, Theorem~\ref{thm:bip3} (3) says that 
$s\le 12$ while the bound ~\eqref{ineq:drt} says that $s\le 13$.
\end{remark}
\section{Concluding remarks}\label{sec:6}
In this paper, we gave upper bounds on the size of transitive subtournaments 
in a digraph. 
In particular, we obtained an analogy of Hoffman's bound to digraphs. 
Furthermore, we partially improved  the Hoffman type bound for doubly regular 
tournaments by applying block intersection polynomials, and thus we found a new application of block intersection polynomials. 

Interesting problems which are worth looking into as future works are listed below. 
\begin{itemize}
\item We could not find any example of digraphs attaining the equality of the upper bound~\eqref{thm:hoff1}. Find examples of such digraphs or  prove the nonexistnce of such digraphs. 
\item  Bachoc et al.~\cite{BMR} improved the Hoffman bound for the Paley graphs by using the properties of quadratic residues of finite fields. Find an analogy of Bachoc et al.'s result to the Paley tournaments. 
\item Hoffman~\cite{H} gave a bound on the chromatic number of a graph. In particular, it was  proved that for any graph $G$,  $\chi(G)\ge 1-\lambda_\max/\lambda_{\min}$, where $\chi(G)$ is the chromatic number of $G$, and $\lambda_{\max}$ and $\lambda_\min$ are the maximum and minimum eigenvalues of $G$, respectively. 
Find an analogy of this bound to digraphs. Here, the chromatic number of a digraph is defined to be the minimum number of 
disjoint transitive subtournaments covering all vertices. 
\end{itemize} 

\end{document}